\documentclass[a4paper,10pt,nointlimits]{article}
\usepackage{amssymb,amsmath,amsthm,verbatim}
\usepackage[dvips,final]{graphicx}
\usepackage{psfrag}
\usepackage{natbib}
\newcommand{\R}{\mathbb{R}}
\newcommand{\Prob}{\mathbb{P}}
\renewcommand{\Pr}{\Prob}

\newcommand{\E}{\mathbb{E}}

\newcommand{\Fc}{\mathcal{F}}

\newcommand{\half}{\frac{1}{2}}

\newtheorem{theorem}{Theorem}
\newtheorem{prop}[theorem]{Proposition}
\newtheorem{lemma}[theorem]{Lemma}
\newtheorem{cor}[theorem]{Corollary}

\theoremstyle{definition}

\newtheorem{remark}[theorem]{Remark}
\newtheorem{definition}[theorem]{Definition}

\newcommand{\indic}[1]{\boldsymbol{1}_{\{\ensuremath{#1}\}}}

\DeclareMathOperator*{\argmax}{argmax}
\DeclareMathOperator*{\argmin}{argmin}
\newcommand{\as}{a.s.}

\newcommand{\eg}{e.g.}

\newcommand{\tdL}{\tilde{L}}
\newcommand{\tdB}{\tilde{B}}

\begin{document}
\title{A unifying class of Skorokhod embeddings: connecting the Az\'{e}ma-Yor and Vallois embeddings}
\author{A.~M.~G.~Cox\thanks{e-mail:
        \texttt{amgc500@york.ac.uk}; web:
        \texttt{www-users.york.ac.uk/$\sim$amgc500/}}\\
        Department of Mathematics,\\
        University of York,\\
        York Y010 5DD, U.~K.        
 \and D.~G.~Hobson\thanks{e-mail:
        \texttt{dgh@maths.bath.ac.uk}; web:
        \texttt{www.bath.ac.uk/$\sim$masdgh/}},\\
        Department of Mathematical Sciences,\\
        University of Bath,\\
        Bath BA2 7AY, U.~K.}
\maketitle

\begin{abstract}
In this paper
we consider the Skorokhod embedding problem in Brownian motion. In 
particular, we give a solution 
based on the 
local time at zero of a variably skewed Brownian motion related to the 
underlying Brownian motion. Special cases of the construction include 
the Az\'{e}ma-Yor and Vallois embeddings. In turn, the 
construction has an interpretation in the Chacon-Walsh framework.

%
%
\end{abstract}

\section{Introduction}
The Skorokhod embedding problem was first proposed, and then
solved, by \citet{Skorokhod:65}, and may be described thus:
\begin{quote}
Given a Brownian motion $(B_t)_{t \ge 0}$ and a centred target law
$\mu$ can we find a stopping time $T$ such that $B_T$ has
distribution $\mu$?
\end{quote}
Skorokhod gave an explicit construction of the stopping time $T$
in terms of independent random variables, and in doing so showed
that any zero-mean probability measure may be embedded in Brownian
motion. Since the problem was posed many more solutions have been
given, see for example \citet{Dubins:68}, \citet{Root:69} and
\citet{Bass:83}, and the comprehensive survey article of 
\citet{Obloj:04}.

In this paper we are concerned with the connection between two of the 
well-known constructions due to \citet{AzemaYor:79} and 
\cite{Vallois:83}, and the relationship between these constructions and 
the construction of \cite{ChaconWalsh:76}. 
The \citet{AzemaYor:79} construction defines a stopping time
\[
T_{AY} = \inf\{t \ge 0 : B_t \le \psi(S_t)\}
\]
where $S_t$ is the maximum of the Brownian motion $B$ and $\psi$ is an 
increasing function dependent on $\mu$. This is an embedding of $\mu$ 
and has a nice optimality property; namely that it maximises the 
distibution of the maximum for embeddings under which $B_{t \wedge T}$ 
is uniformly integrable (see \citet{AzemaYor:79b}). The other well-known construction of interest to us here is the 
embedding of \citet{Vallois:83}. This construction uses the local 
time $L^B_t$ of the Brownian motion at 0. In particular, for certain 
functions $\alpha,\beta$, 
\[
T_V = \inf\{t \ge 0 : B_t \not\in (\beta(L^B_t),\alpha(L^B_t))\}
\]
is an embedding. In \citet{Vallois:92} it is shown that for a convex 
function $\Psi$ a specific 
choice of $\alpha$ and $\beta$ 
minimises $\E \Psi(L^B_T)$ among the class of all embeddings of $\mu$. Moreover the 
same choice of $\alpha$ and $\beta$ is optimal for each $\Psi$.
It is this specific choice of $\alpha$ and
$\beta$, for which
$\alpha$ is decreasing and $\beta$ increasing,
which interests us here.

It is not surprising given L\'evy's result on the equivalence in law between 
$(S_t,S_t-B_t)$ and $(L^B_t,|B_t|)$ that there is a strong connection 
between the two constructions (not least when $\mu$ is symmetric); 
however the aim of this paper is to demonstrate that both embeddings can 
be seen as examples of a much broader class of embeddings. In 
particular, our embeddings will be characterised by a function $G: \R^+ 
\to \R$, and related to variably skewed Brownian motion. Skewed Brownian 
motion \citep{HarrisonShepp:81} can be thought of as a Brownian motion 
with excursions being positive with probability $p$ and negative with 
probability $1-p$; for a variably skewed Brownian motion 
\citep{Barlowetal:00} the probability of positive excursions will vary 
with the local time at zero of the process.

Given a Brownian motion $B_t$, a variably skewed Brownian motion $X_t$ 
is the solution to
\begin{equation} \label{eqn:Xdefn1}
X_t = B_t - G(L_t^X)
\end{equation}
where $L_t^X$ is the local time of $X_t$ at 0 and $G$ is a function with Lipschitz constant at most $1$. Our embedding will then 
be a stopping time of the form
\begin{equation} \label{eqn:Tdefn1}
T = \inf\{t \ge 0 : X_t \not\in (b(L_t^X),a(L_t^X))\} 
= \inf\{t \ge 0 : B_t \not\in (\beta(L_t^X),\alpha(L_t^X))\}
\end{equation}
for suitable pairs of functions $(a,b)$ and $(\alpha,\beta)$. When $G \equiv 0$, $X_t$ is 
just $B_t$ (this remark can be restated in a more complicated fashion by 
saying that $X$ is a skew-Brownian motion with constant skewness 
parameter 1/2) and the stopping time is the same as Vallois' stopping 
time. 
When $G(x) = x$ for all $x$, we find that $X_t$ has only negative 
excursions and $L_t^X = \sup_{s \le t} B_s$; in particular, the 
embedding is now of the form
\[
T = \inf\{t \ge 0 : B_t \le  \psi(S_t))\}
\]
for some $\psi$ and we recover the Az\'{e}ma-Yor stopping time.

The embeddings we construct have an interpretation within the balayage 
construction of \cite{ChaconWalsh:76}. As shown in \citet{Meilijson:83}, 
the Az\'{e}ma-Yor construction is already known to fit in this framework. 
We shall show how our construction can be interpreted in a 
similar way.

The paper proceeds in the following way: we shall begin by 
establishing the important quantities in the construction, and some 
relationships between them; then we show in Theorem~\ref{thm:varskew} 
that there is a strong solution to \eqref{eqn:Xdefn1}, and hence that 
given the Brownian motion, we may indeed construct $X$. The main result 
of this paper is Theorem~\ref{thm:embeds}, where we use excursion 
arguments to confirm that the stopping time is an embedding, and that 
$B_{t \wedge T}$ is uniformly integrable. In Section~\ref{sec:examples} we discuss the connection with the Az\'{e}ma-Yor and Vallois embeddings, and describe a new embedding in the general framework which maximises $\E\Psi(L_T^x)$ for $\Psi$ convex and $L_T^x$ the local time of the Brownian motion at a non-zero level $x$.

\section{The Construction}

Our aim is to show that we can construct embeddings of a target 
distribution $\mu$ in a Brownian motion $(B_t)_{t \ge 0}$ with $B_0 = 0$, 
where the embedding is determined by a curve $\gamma:[0,\zeta] \to \R 
\times \R_+$, $\gamma(s) = (F(s),h(s))$ and what is in some sense the 
local time of the Brownian motion on the function $F$. The 
parameterisation of the curve is essentially arbitrary, but for 
simplicity we might choose $s$ to be the arc length. The construction we 
provide will have as special cases the embeddings of \citet{AzemaYor:79} 
and \citet{Vallois:83}. In particular, we will see that the Az\'{e}ma-Yor 
construction corresponds to the choice $F(s) = h(s)$ and the 
Vallois construction to the choice $F \equiv 0$.

In order to perform the construction, we need to make some assumptions 
on the behaviour of $\gamma$. The exact form of the conditions is 
dependent on the target distribution, and further notation is required 
to state them clearly, however the assumptions have a 
straight-forward graphical interpretation.

As with many embeddings, there are strong connections with potential 
theory. 
We define the function
\begin{equation} \label{eqn:cdefn}
c(x) = \int_\R |y-x| \, \mu(dy) = \E^\mu |X-x|.
\end{equation}
This is related to a common definition of the potential 
of a Brownian motion (in 1-dimension) via 
$c(x) = -u_{\mu}(x)$ where $u_{\mu}$ is the potential of $\mu$ (see \eg{} \citet{Chacon:77}).

\begin{prop}\label{prop:cprop}
Let $\mu$ be an integrable, centred target distribution. 
Then the function $c$ as defined in \eqref{eqn:cdefn} 
has the following useful properties:
\begin{enumerate}
\item
$c$ is a  positive, continuous, convex function such that
\(
c(x) \ge |x|;
\)
\item
as $x \to \pm \infty$, $c(x) - |x| \to 0$;
\item
$c$ has left and right derivatives, given by:
\begin{equation*}
c_-'(x) = 2 \mu((-\infty,x)) -1, \quad c_+'(x) = 2 \mu((-\infty,x]) - 1,
\end{equation*}
so in particular, $c$ is differentiable Lebesgue-almost everywhere. 
The set $\{x : \mu(\{x\})>0\}$ is precisely the set where $c$ is not 
differentiable.
\end{enumerate}
\end{prop}

We now discuss the set of admissible curves $\gamma$.
We suppose that the functions $F$ and $h$ are both absolutely 
continuous, and hence have derivatives Lebesgue-almost everywhere. We 
shall assume for convenience that the parameterisation is such that $s$ 
is the arc-length, however this is not necessary, and the extension to 
an arbitrary parameterisation is immediate provided that 
$Leb(\{s : h'(s) = F'(s) = 0\}) = 0$. We 
also define $\zeta = \sup\{s>0 : s \le c(F(s))\}$, and require that the 
curve $\gamma(s) = (F(s),h(s))$ 
hits c, at least in the limit $s \to 
\zeta$ (although this could occur at $\pm \infty$).

\psfrag{phi}{$\phi(s)$}
\psfrag{y}{$h(s)$}
\psfrag{theta}{$\theta(s)$}
\psfrag{F}{$F(s)$}
\psfrag{c}{$c(x)$}
\begin{figure}[t]
\begin{center}
\includegraphics[width=\textwidth,height=3in]{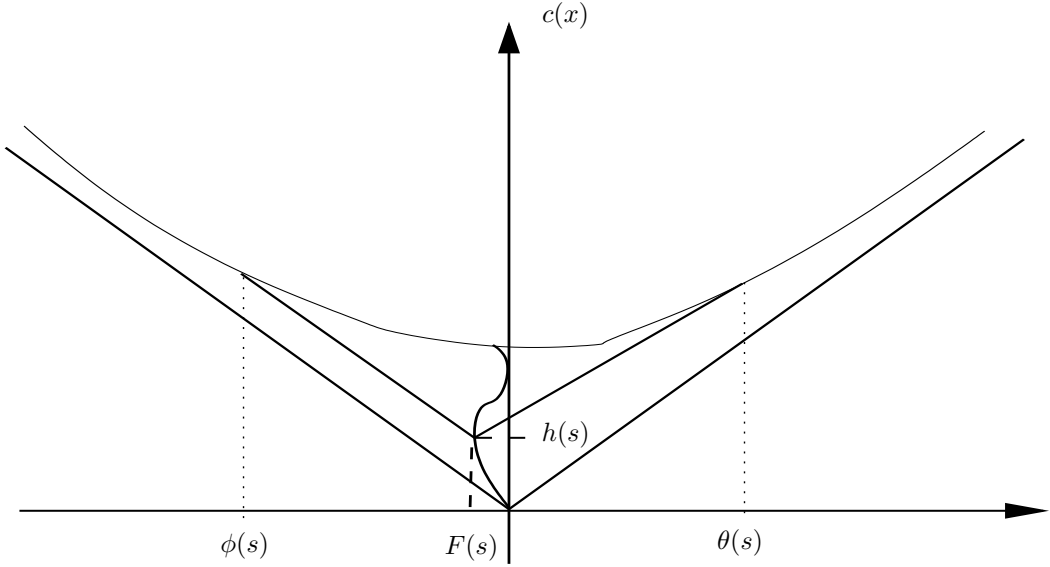}
\caption{\label{fig:basic} The `potential' $c(x)$ for a distribution 
$\mu$. The embedding is determined by the curve $\gamma = (F,h)$
and $\theta(s)$ and $\phi(s)$ are the $x$-coordinates of the 
points at which the tangents to $c$ intersect $\gamma$ at $\gamma(s)$.
} \end{center}
\end{figure}

The curve $\gamma$ is key to defining the other primary variables.
The functions $\theta, \phi$ are defined in terms of $c$, $F$ and $h$ 
so that for $s \ge 0$,
\begin{eqnarray}
\theta(s) & = & \sup\left\{x \in \R:c(\lambda) \ge h(s) + (\lambda - F(s)) 
\frac{c(x)-h(s)}{x-F(s)} \quad \forall \lambda \le x\right\} 
\label{eqn:thetadefn},\\
\phi(s) & = & \inf\left\{x \in \R: c(\lambda) \ge h(s) + (\lambda- F(s)) 
\frac{h(s) - c(x)}{F(s) - x} \quad \forall \lambda \ge x\right\}. 
\label{eqn:phidefn}
\end{eqnarray}
Alternative definitions of $\theta$ and $\phi$ are given by
\begin{eqnarray}
\theta(s) & = & 
\argmin_{x>F(s)}\left\{\frac{c(x) - h(s)}{x-F(s)}\right\},
\label{eqn:thetadefn2}\\
\phi(s) & = & 
\argmax_{x < F(s)}\left\{\frac{h(s)-c(x)}{F(s) - x}\right\}
\label{eqn:phidefn2},
\end{eqnarray}
but note that the $\argmin$ and $\argmax$ in \eqref{eqn:thetadefn2} and 
\eqref{eqn:phidefn2} may not be uniquely defined. There are no such 
problems with the definitions in \eqref{eqn:thetadefn} and 
\eqref{eqn:phidefn}. 

In order for the construction to make sense, we must restrict $\gamma$
so that for any particular value of $s$ the future of the curve will lie
inside the (closed) area constructed via the tangents at
$\theta(s)$ and $\phi(s)$ and the potential $c$. We note that under this
restriction, the curve can have $h'(s) < 0$, which is why we cannot
paramatrise the curve by its height. Our first aim is to express this
idea mathematically. Define 
\begin{eqnarray}
R(s) & =  & \frac{c(\theta(s)) - h(s)}{\theta(s) - F(s)} =
\inf_{x>F(s)} \frac{c(x)-h(s)}{x-F(s)},\label{eqn:Rdefn}\\
S(s) & = & \frac{h(s)-c(\phi(s))}{F(s) - \phi(s)} =
\sup_{x<F(s)} \frac{h(s)-c(x)}{F(s) - x},\label{eqn:Sdefn}\\
\Gamma(s) & = & \frac{R(s) - S(s)}{2}.\label{eqn:Gammadefn}
\end{eqnarray}
In particular, $R(s)$ and $S(s)$ are the gradients of tangents 
to $c$ which hit $\gamma=(F,h)$ at $(F(s),h(s))$. As well as $\Gamma$ 
representing half the difference in the slopes it will also have an 
interpretation as the probability that the embedding has not yet 
occured. If we use the notation 
$f_+'$ to denote the right derivative of $f$ then the condition
on $\gamma$ is
\begin{equation}
h_+'(s) \ge \max(F_+'(s) S(s),F_+'(s) R(s)). \label{eqn:hdashcond}
\end{equation}

\psfrag{phi}{$\phi(s)$}
\psfrag{y}{$h(s)$}
\psfrag{theta}{$\theta(s)$}
\psfrag{F}{$F(s)$}
\psfrag{c}{$c(x)$}
\psfrag{B}{$B_t$}
\psfrag{X}{$X_t$}
\psfrag{Y}{$Y_t$}
\psfrag{L}{$L_t^X$}
\psfrag{Ytis}{$L_t^X = H(Y_t)$}
\psfrag{Xtis}{$X_t = B_t - G(L_t^X)$}
\begin{figure}[t]
\begin{center}
\includegraphics[width=\textwidth,height=2in]{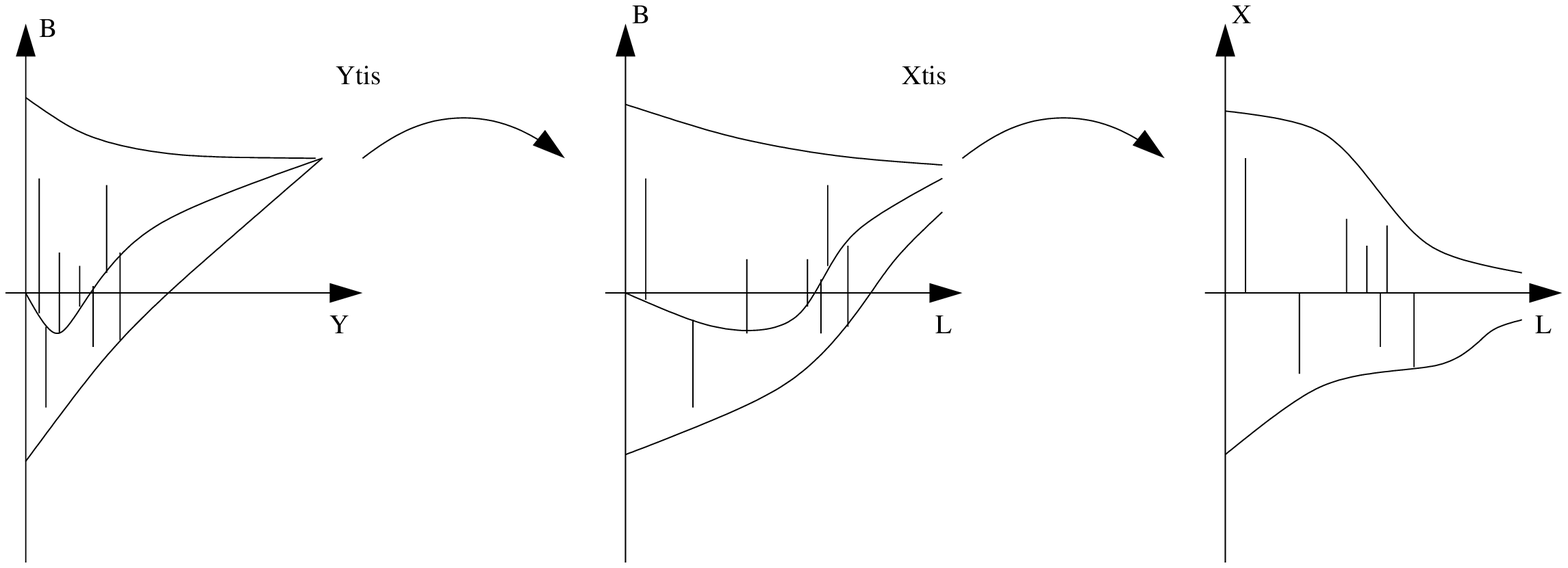}
\caption{\label{fig:change} Three different interpretations of the 
stopping time; the first graph shows the functions $\theta(s), \phi(s)$ 
and $F(s)$; the horizontal axis represents the current `position along 
the curve' or the `local time' of $B$ on $F$ labelled by $Y_t$ --- our 
primary aim is to understand how to interpret this. In the second graph 
we transform the horizontal scale using a function $H$ to be determined, 
so that the horizontal axis is the local time $L_t^X$. For the final 
graph we transform the vertical scale by subtracting $G=F \circ H^{-1}$, 
so that the vertical axis is $X_t = B_t - G(L_t^X)$.  The 
vertical lines are the excursions of $B$ from $F$, or excursions of $X$ 
from zero. In each interpretation the stopping time can be though of as 
the first time an 
excursion goes outside the region bounded by the upper and lower 
curves.}
\end{center}
\end{figure}

The construction of the stopping time can be considered via the 
transformations shown in Figure~\ref{fig:change}. We want to be 
able 
to make sense of the `local time' of $B_t$ on $F$, which we denote by 
$Y_t$. For the present we are not specific about the process $Y$ except 
that it is an increasing process which only grows when `$B$ is on $F$' 
in a sense which will become clear. The full definition of $Y$ can be 
obtained from reversing the construction process we are about to 
describe.

The curves in the left graph of Figure~\ref{fig:change} (starting with the highest) represent $\theta, F$ and $\phi$. In order to define the candidate stopping time we need to transform the 
horizontal and vertical scales of this graph. Firstly, we 
map the horizontal axis by $H^{-1}$ (where $H$ will be defined later), 
so that the curves in the centre graph are plots of $F(H^{-1}(l)), \theta(H^{-1}(l)), 
\phi(H^{-1}(l)) \equiv (G(l), \alpha(l), \beta(l))$. Secondly, we 
transform the vertical scale by subtracting $G(l)$ from the vertical co-ordinate, so that the third graph shows $a(l)= \alpha(l) - G(l)$ and 
$b(l)= \beta(l)-G(l)$. Also shown in Figure \ref{fig:change} (by the 
vertical lines) are the excusions of the Brownian motion away from $F$. In the first and second graph these excursions are plotted in the co-ordinates $(Y_t,B_t)$ and $(L_t^X,B_t)$ respectively.
In the final plot we have a process $X_t$ making excursions from $0$, 
and it now makes sense to define the horizontal axis to be the 
(symmetric) local time of $X_t$ at 0. In particular, if the process in 
the final graph has local time $L_t^X$, then $X$ must satisfy
\[
X_t = B_t - G(L_t^X).
\]
We will later show (Theorem~\ref{thm:varskew}) that this equation has a 
strong solution. The candidate embedding is then of the form
\[
T = \inf \{ t \ge 0 : X_t \not \in (b(L_t^X),a(L_t^X))\}
\]
where $a(l) = \theta(H^{-1}(l)) - F(H^{-1}(l))$ and $b(l) = 
\phi(H^{-1}(l)) - F(H^{-1}(l))$. This stopping rule can then be re-expressed in the form \eqref{eqn:Tdefn1}.

Having described the form and philosophy of the construction it remains 
to determine the transformation function $H$, and then to show that the 
algorithm outlined above leads to a genuine embedding of the target law 
$\mu$. The advantage of representing the construction via the diagrams 
in Figures~\ref{fig:basic} and \ref{fig:change} is to show how these 
representations are special cases of the \citet{ChaconWalsh:76} 
embedding. 

The conditions on the curve $\gamma$ ensure that the 
function $\theta$ is 
non-increasing, and $\phi$ is non-decreasing. It also follows from 
the definitions that $\theta$ and $\phi$  
are both monotonic and left-continuous. Note that $\alpha$ and $\beta$ 
share these properties, but that $a$ and $b$ need not be monotonic. 

The first results discuss the properties of the functions $R$ and $S$ 
introduced in \eqref{eqn:Rdefn} and \eqref{eqn:Sdefn}.
Note that (by  Proposition~\ref{prop:cprop}(iii)) we have the inequalities
\begin{eqnarray}
1-2\mu([\theta(s),\infty)) & \le R(s) \le
& 1 - 2 \mu((\theta(s),\infty)),\label{eqn:Rbound}\\
2\mu((-\infty,\phi(s))) - 1 & \le S(s) \le & 
2 \mu((-\infty,\phi(s)]) -1.\label{eqn:Sbound}
\end{eqnarray}

\begin{lemma} \label{lem:A}
The functions $R$ and $S$ are absolutely continuous on closed subsets 
of $[0,\zeta)$
and
\begin{eqnarray*}
R(s) & = & \int_0^s \frac{-1}{\theta(z) - F(z)} 
\left[ h'(s) - F'(z) R(z)\right]\, dz,\\
S(s) & = & \int_0^s \frac{1}{F(z) - \phi(z)} 
\left[h'(s)-F'(z) S(z) \right]\, dz .
\end{eqnarray*}
\end{lemma}
\begin{proof} The function
$R(s)$ is the gradient of the tangent to $c$ that hits 
$\gamma$ at $(F(s),h(s))$. It is clear that $R$ is 
non-increasing 
(by the 
assumptions on $\gamma$), and for every gradient there is a unique 
tangent which 
hits $\gamma$ at some point; hence $R$ must be continuous.

Also, for $\delta$ small (we insist that $F(s - \delta) < \theta(s)$) we 
have $F(s-\delta) \le F(s) + A\delta$ and 
$h(s-\delta) \ge h(s) - B\delta$ for some $A,B>0$, and by \eqref{eqn:Rdefn}
\begin{eqnarray*}
R(s - \delta) & \le & \frac{c(\theta(s)) - h(s) 
+ B\delta}{\theta(s) - F(s) - A\delta}\\
& = & \frac{c(\theta(s)) - h(s)}{\theta(s) - F(s)} 
+ \frac{c(\theta(s)) - h(s)}{\theta(s) - F(s)-A\delta} 
- \frac{c(\theta(s)) - h(s)}{\theta(s) - F(s)} 
+ \frac{B\delta}{\theta(s) - F(s)-A\delta}\\
& = & R(s) + \left[ A\frac{c(\theta(s)) - h(s)}{\theta(s) - F(s)} + 
B\right] \frac{\delta}{\theta(s) - F(s) - B\delta}
\end{eqnarray*} 
so that $R$ is locally Lipschitz on $[0,\zeta)$ and absolutely 
continuous on $[0,z]$ for $z < \zeta$. In particular it is 
differentiable Lebesgue-almost everywhere, and if $R'$ is a version of 
the derivative of $R$, then we must have
\[
R(s) = 1 + \int_0^s R'(z) \, dz
\]
for $s<\zeta$ (and hence by continuity, for $s = \zeta$).

We can calculate the left-derivative of $R(s)$, using the fact that 
$\theta$ is left continuous and decreasing, to be:
\begin{equation*}
\frac{d^-}{ds}R(s) = -\frac{h_-'(s)}{\theta(s) - F(s)} + 
\frac{F_-'(s)R(s)}{\theta(s) - F(s)} + \frac{\theta_-'(s)}{\theta(s) - F(s)} 
\left[ c_+'(\theta(s)) 
- \frac{c(\theta(s)) - h(s)}{\theta(s) - F(s)} \right].
\end{equation*}
However the second term will disappear: if there is no atom at 
$\theta(s)$, it is easy to see that
\begin{equation} \label{eqn:cplusdiff}
c_+'(\theta(s)) = \frac{c(\theta(s)) - h(s)}{\theta(s) - F(s)}
\end{equation}
while if there is an atom at $\theta(s)$, $\theta_-'(s) = 0$ (unless $s$ 
is the largest $s$ corresponding to the atom, however then 
\eqref{eqn:cplusdiff} still holds).

Consequently
\[
\frac{d^-}{ds}R(s) 
= -\frac{h_-'(s)}{\theta(s) - F(s)} + \frac{F_-'(s)R(s)}{\theta(s) - F(s)}
\]
and the result holds. We can deduce the statement for $S(s)$ similarly.
\end{proof}

\begin{cor}
The function $\Gamma$ is absolutely continuous on closed subsets of 
$[0,\zeta)$ and Lebesgue-almost everywhere differentiable with
\begin{equation*}
\frac{d\Gamma}{ds} = \frac{-1}{2(\theta(s) - F(s))}(h'(s)-F'(s) R(s)) 
- \frac{1}{2(F(s) - \phi(s))} (h'(s)-F'(s) S(s))
\end{equation*}
when the derivative exists.
\end{cor}

\begin{prop}
Define
\begin{eqnarray} \label{eqn:hdefn}
H(s) & = & \int_0^s \frac{h'(z)}{\Gamma(z)}\, dz - \int_0^s 
F'(z)\frac{R(z) + S(z)}{R(z)-S(z)}\, dz\\ 
& = & \int_0^s \left[2h'(z) - 
F'(z)\left(R(z) + S(z)\right)\right]\frac{dz}{R(s) - S(s)}.
\end{eqnarray}
Then $H$ is strictly increasing and continuous on $[0,\zeta)$, with 
continuous inverse $H^{-1}$.
\label{prop:H}
\end{prop}

\begin{proof}
This follows from \eqref{eqn:hdashcond} and the fact that 
$R(z)=S(z)$ only if $z=\zeta$.
\end{proof}

\begin{prop} \label{prop:Hprop}
The function $G(l) = F \circ H^{-1}(l)$ 
is absolutely continuous with $|G'(l)| \le 1$.
\end{prop}
\begin{proof}
Both $F$ and $H^{-1}$ are differentiable almost everywhere, so that 
\[
G'(l) = \frac{F'(H^{-1}(l))}{H'(H^{-1}(l))}.
\]
It is sufficient to show that $|G'(l)| \le 1$, or equivalently that $-H' 
\le F' \le H'$. However:
\begin{equation*}
H' - F' = \frac{2h'}{R-S} - \frac{2F'R}{R-S} = \frac{2}{R-S} (h'-F'R) \ge 0,
\end{equation*}
and
\begin{equation*}
H' + F' = \frac{2h'}{R-S} - \frac{2F'S}{R-S} = \frac{2}{R-S} (h'-F'S) \ge 0,
\end{equation*}
by \eqref{eqn:hdashcond}.
\end{proof}

The following lemma follows immediately from 
Lemma~\ref{lem:A} and Proposition~\ref{prop:H}.

\begin{lemma} \label{lem:F}
We have the relationships
\[
\frac{H'-F'}{\theta-F} = -\frac{R'}{\Gamma}, 
\quad \frac{H'+F'}{F-\phi} = \frac{S'}{\Gamma}.
\]
\end{lemma}


\section{The verification.}
We have now completed the definitions of the auxilliary functions. It 
remains to show, firstly that the process given by \eqref{eqn:Xdefn1} is 
well-defined, and secondly that the stopping time in \eqref{eqn:Tdefn1}
embeds the law $\mu$.

\begin{theorem}\label{thm:varskew}
Consider the equation
\begin{equation} \label{eqn:Xdefn}
X_t = B_t - G(L_t^X)
\end{equation}
where $B_t$ is a Brownian motion, and $L_t^X$ is the local time at 0 of 
the process $X$. Then there exists a strong solution to 
\eqref{eqn:Xdefn}. Further, the solution is a variably skewed Brownian 
motion, with probability of positive excursions $p$ given by
\[
p(l) = \frac{1 - G'(l)}{2}
\]
when the local time $L_t^X = l$.
\end{theorem}

The existence of a strong solution is shown in \citet{Barlowetal:00} 
under the slightly stronger condition that $|G'(l)| < 1$. They prove 
that the solutions of the equation are pathwise unique and that there 
exists a weak solution, which, by the Yamada-Watanabe Theorem, implies 
the existence of a strong solution. (The form of the Yamada-Watanabe 
Theorem given in \eg{} \citet{RevuzYor:99}[IX 1.7] is stated for an 
ordinary SDE, but as Barlow et al note, can be extended to cover the 
case of interest.) The proof of Theorem~3.1 in \citet{Barlowetal:00} 
carries over without alteration to the case where $|G'(l)| \le 1$, so it 
only remains to show the existence of a weak solution in the more 
general case. The statement about the probability of excursions follows 
on consideration of the case where $G(x) = \beta x$ for $|\beta| \le 1$.

Consequently, Theorem~\ref{thm:varskew} will follow from the following 
result:
\begin{lemma} \label{lem:weaksoln}
There exists a weak solution to \eqref{eqn:Xdefn}.
\end{lemma}

The proof of this lemma will follow the proof of Proposition~2.2 of 
\citet{Barlowetal:00}; the major difference is that it is necessary to 
use two functions $n,m$ instead of the single function $\phi$ used by 
Barlow et al.

\begin{proof}
For $x \geq 0$ define $v(x)$ to be the solution of
\[
\frac{dv}{dx} = \left( \frac{2}{1 + (G'(x))^2}\right)^\half
\]
with $v(0)=0$, so that $1 \le v'(x) \le \sqrt{2}.$ Define $y$ to be 
the 
inverse of $v$. Then $\frac{1}{\sqrt{2}} \le y'(x) \le 1$. Define also
\begin{eqnarray}
n(x) & = & \frac{y(x) - G(y(x))}{2}\label{eqn:ndefn};\\
m(x) & = & \frac{y(x) + G(y(x))}{2}\label{eqn:mdefn}.
\end{eqnarray}
In particular, $0 \le n' \le y' \le 1$ and $0 \le m' \le y' \le 1$.

Now define
\begin{equation*} 
R(x,l) = \begin{cases} 2 n'(l)x &:x \ge 0\\
2m'(l)x &:x < 0
\end{cases}.
\end{equation*}
Let $\tdB_t,\tdL_t$ be a Brownian motion and its local time. 
By applying \citet{RevuzYor:99}[VI.4.3] separately to 
$R_+(\tdB_t,\tdL_t) = 2n'(\tdL_t)(\tdB_t\vee 0)$ and 
$R_-(\tdB_t,\tdL_t) = m'(\tdL_t)(\tdB_t \wedge 0)$ we deduce
\begin{eqnarray*}
R(\tdB_t,\tdL_t) & = & 2 \int_0^t(n'(\tdL_s)\indic{\tdB_s \ge 0} 
+ m'(\tdL_s)\indic{\tdB_s < 0}) \, d\tdB_s 
+ \int_0^t (n'(\tdL_s)-m'(\tdL_s))\, d\tdL_s\\
& = & 2 \int_0^t(n'(\tdL_s)\indic{\tdB_s \ge 0} 
+ m'(\tdL_s)\indic{\tdB_s < 0}) \, d\tdB_s + n(\tdL_t)-m(\tdL_t),\\
|R(\tdB_t,\tdL_t)| & = & 2 \int_0^t(n'(\tdL_s)\indic{\tdB_s \ge 0} 
- m'(\tdL_s)\indic{\tdB_s < 0}) \, d\tdB_s + \int_0^t (n'(\tdL_s)
+m'(\tdL_s))\, d\tdL_s\\
& = & 2 \int_0^t(n'(\tdL_s)\indic{\tdB_s \ge 0} 
- m'(\tdL_s)\indic{\tdB_s < 0}) \, d\tdB_s + n(\tdL_t)+m(\tdL_t).
\end{eqnarray*}
The stochastic integrals in both of the above equations have 
quadratic variation
\[
A_t = 4 \int_0^t(n'(\tdL_s)^2\indic{\tdB_s \ge 0} 
+ m'(\tdL_s)^2\indic{\tdB_s < 0}) \, ds.
\]

Since $n'$ and $m'$ are bounded, $A_t < \infty$ \as{}, and because $n' 
\vee 
m' \ge \half y' \ge \frac{1}{2\sqrt{2}}$ we also have $A_t \uparrow 
\infty$ \as{}.
Consequently, defining $T_t = \inf\{s : A_s >t\}$, we deduce that
\begin{eqnarray*}
W_t & = & 2 \int_0^{T_t}(n'(\tdL_s)\indic{\tdB_s \ge 0} 
+ m'(\tdL_s)\indic{\tdB_s < 0}) \, d\tdB_s,\\
\tilde{W}_t & = & 2 \int_0^{T_t}(n'(\tdL_s)\indic{\tdB_s \ge 0} 
- m'(\tdL_s)\indic{\tdB_s < 0}) \, d\tdB_s
\end{eqnarray*}
are Brownian motions.

Since $L_t = n(\tdL_{T_t}) + m(\tdL_{T_t})$ is increasing only on the 
set $\{t:R(\tdB_{T_t},\tdL_{T_t}) = 0\}$ and $|R(\tdB_{T_t},\tdL_{T_t})| 
= \tilde{W}_t + L_t$, $L_t$ must be the symmetric local time at 0 of 
$R(\tdB_{T_t},\tdL_{T_t})$.

However, by \eqref{eqn:ndefn} and \eqref{eqn:mdefn}, we have
\[
R(\tdB_{T_t},\tdL_{T_t}) = W_t - G(L_t),
\]
and hence a weak solution to \eqref{eqn:Xdefn}.
\end{proof}

The important fact is that $X_t$ is a strong solution, so that it 
is adapted to the filtration $\Fc_t^B$ of the Brownian motion --- in 
particular, given the Brownian motion, it is possible to construct $X_t$ 
and its local time $L_t^X$, even if we do not have an explicit 
construction (see however Remark~3.5 of Barlow et al.).

Given a curve $\gamma$, we can define $\theta$ and $\phi$ as in 
\eqref{eqn:thetadefn} and \eqref{eqn:phidefn}, and thence define $R$, 
$S$, 
$H$ and $G= F 
\circ H^{-1}$. Define also
\begin{eqnarray}
a(l) & = & \theta(H^{-1}(l)) - F(H^{-1}(l)) \label{eqn:adefn} \\
b(l) & = & \phi(H^{-1}(l)) - F(H^{-1}(l)) \label{eqn:bdefn}
\end{eqnarray}
and let $X$ be given by the solution to \eqref{eqn:Xdefn}.

\begin{definition}
Define the stopping time
\begin{equation}\label{eqn:Tdefn}
T = T_{\gamma} = \inf\{t \ge 0: X_t \not\in (b(L_t^X),a(L_t^X))\}.
\end{equation}
$T$ will be our embedding of $\mu$.
\end{definition}

The key idea will be that we can use excursion theory to calculate the 
distribution of the stopped process $B_T$. 
We shall first need the following result:

\begin{prop}
Let $Y_t = H^{-1}(L^X_t)$. The the law of $Y_T$ is given by
\begin{equation}\label{eqn:Yprob} 
\Pr(Y_T \ge s) = \Gamma(s).
\end{equation}

\end{prop}
\begin{proof}
By an excursion theory argument, the rate $r(l)$ at which positive 
excursions 
from zero hit $a(l)$ is the rate at which an excursion from 0 has 
maximum modulus 
$a(l)$ multiplied by the probability the excursion is positive:
\[ 
r(l) = \frac{1}{a(l)} \times \frac{1-G'(l)}{2} = \frac{1 - 
G'(l)}{2(\theta(H^{-1}(l)) - F(H^{-1}(l)))}.
\]
Similarly, the rate of negative excursions is
\[
\frac{1}{|b(l)|} \times \frac{G'(l)+1}{2} = 
\frac{G'(l)+1}{2(F(H^{-1}(l)) - \phi(H^{-1}(l)))}.
\]
Consequently
\begin{eqnarray*}
\Pr(Y_T \ge s) & = & \Pr(L_T^X \ge H(s))\\
& = & \exp\left\{-\int_0^{H(s)}\left[
\frac{1 - G'(l)}{2(\theta(H^{-1}(l)) - F(H^{-1}(l)))} 
+ \frac{G'(l)+1}{2(F(H^{-1}(l)) - \phi(H^{-1}(l)))}\right] \, dl\right\}\\
& = & \exp\left\{-\int_0^s \half\left[ 
\frac{H'(z) - F'(z)}{\theta(z) - F(z)} 
+ \frac{F'(z) + H'(z)}{F(z) - \phi(z)}\right] \, dz\right\}\\
& = & \exp\left\{\int_0^s \frac{\Gamma'(z)}{\Gamma(z)} \, dz\right\}\\
& = & \Gamma(s)\\
\end{eqnarray*}
where the penultimate equality is a consequence of Lemma~\ref{lem:F}.
\end{proof}

\begin{theorem} \label{thm:embeds}
The stopping time $T$ is a UI-embedding: 
that is $B_T \sim \mu$ and the process $B_{t \wedge T}$ 
is uniformly integrable.
\end{theorem}

\begin{remark}
We note that $T$ being a UI-embedding is equivalent to $T$ being minimal
in the sense of \citet{Monroe:72}; that is $T$ is minimal if whenever 
$S$ is a stopping time with $S \le T$ and $B_S \sim \mu$ we have $S = T$ 
\as{}.  For further discussions on minimality, see 
\citet{CoxHobson:03,Cox:05}.
\end{remark}

\begin{proof}
We first show that $T$ is an embedding. Let $x$ be a point such that 
$x>F(\zeta)$. By definition,
\begin{equation*}
\Pr(B_T \ge x) = \Pr(X_T + G(L_T^X) \ge x).
\end{equation*}
We will only stop at or above $x$ if we have an excursion of $X_t$ above 
$a(L_t^X)$ before the local time is too great; in particular, before 
$\theta(H^{-1}(l)) < x$. When the local time $L_t^X = l$, the rate of 
excursions of $X$ which hit $a(l)$ is $p(l) a(l)^{-1}$. Hence, using 
\eqref{eqn:adefn}, \eqref{eqn:adefn} and \eqref{eqn:Yprob},
\begin{eqnarray*}
\Pr(B_T \ge x) & = & \int_{\{l:\theta(H^{-1}(l))\ge x\}} 
\Pr(L_T \ge l) \frac{p(l)}{a(l)} \, dl\\
& = & \int_{\{z:\theta(z) \ge x\}} \Pr(Y_T \ge z) 
\frac{\half\left(1-\frac{F'(z)}{H'(z)}\right)}{\theta(z) - F(z)}H'(z) \, dz\\
& = & \int_{\{z:\theta(z) \ge x\}} \frac{H'(z)-F'(z)}{2(\theta(z)-F(z))} 
\Gamma(z)\, dz\\
& = & \int_{\{z:\theta(z) \ge x\}} -\frac{R'(z)}{2}\, dz .
\end{eqnarray*}
Since $R$ is continuous and decreasing, it will attain the lower 
bound in \eqref{eqn:Rbound} on the set $\{z:\theta(z) \ge x\}$, 
so that we deduce (since $R(0) = 1$)
\begin{eqnarray*}
\Pr(B_T \ge x) & = & \frac{1 - (1-2\mu([x,\infty)))}{2}\\
& = & \mu([x,\infty)).
\end{eqnarray*}
A similar calculation shows that for $x < \zeta$, 
\[
\Pr(B_T \le x) = \mu((-\infty,x])
\]
and therefore that $B_T \sim \mu$.

To demonstrate uniform integrability, we show that $x \Pr(T > H_x) \to 
0$ as $|x| \to \infty$, where $H_x = \inf\{t \ge 0: B_t = x\}$ (see 
\citet{AzemaGundyYor:79} or \citet{Cox:05}). We consider the case $x 
>0$. If either $F(s)=h(s)$ or $F(s)=-h(s)$ then we are reduced to the 
case of the \citet{AzemaYor:79} embedding (see Remark~\ref{rem:AY}) and 
uniform integrabity follows from standard results. Otherwise, 
the curve $\gamma$ enters the interior of the region between $|x|$ and 
$c(x)$ in Figure~\ref{fig:basic} and then $\sup_{ s \leq \zeta} |F(s)| 
\leq K$ 
for some $K$. Then, for $x>2K$,
\begin{equation*}
\Pr(T>H_x) =  \int_{\{l:\theta(H^{-1}(l)) \ge x\}} \Pr(L_T \ge l) 
\frac{p(l)}{x - F(H^{-1}(l))} \, dl.
\end{equation*}
But $\Pr(L_T \ge 1) \le 1$, $p(l) \le 1$ and $[x - F(H^{-1}(l))]^{-1} < 
(x-K)^{-1} < 2/x$ and so
\[ x \Pr(T>H_x) \leq 2 \mbox{Leb} \{ \{l:\theta(H^{-1}(l)) \ge x \} = 2 H 
(\theta^{-1}(x)) , \]
and as $x \uparrow \infty$, $\theta^{-1}$ tends to zero, and, since $H$ 
is continuous $x \Pr(T>H_x) \rightarrow0$.
\end{proof}

\section{Examples} \label{sec:examples}

\begin{remark} \label{rem:vallois}
The choice $F \equiv 0$ and $h$ increasing can easily be seen to reduce 
to the construction of \citet{Vallois:92} with decreasing functions. 
(\citet{Vallois:92} also has a construction with increasing functions; 
this is not a special case of the current construction.)
\end{remark}

\begin{remark} \label{rem:AY}
Taking $F(s) = h(s) = s/\sqrt{2}$ (or indeed any increasing function of 
$s$), we 
see that $R(s)=1$ and so $H(s) = s$. Consequently we have $G(l) = l$. 
\eqref{eqn:Xdefn} reduces to
\[
X_t = B_t - L_t^X
\]
which is known to have solution $L_t^X = \sup_{s \le t} B_s$, and $X_t$ 
is a reflecting Brownian motion on $(-\infty,0]$. In particular for some 
functions $\phi_1,\phi_2$, our stopping rule is now of the form
\[
T = \inf\{t \ge 0 :X_t \le \phi_1(L_t^X)\} = \inf \left\{t \ge 0 : B_t 
\le \phi_2\left( \sup_{s \le t} B_s\right) \right\},
\]
the second form being notable as the same form as the Az\'{e}ma-Yor 
stopping 
time; it is easy to check the functions are identical.
\end{remark}

We finish by demonstrating how new embeddings can be constructed using the previous results, often with useful optimality properties. The embedding we construct will maximise the expected value of a convex function of the local time at a level $x \neq 0$. (The case where $x=0$ is the Vallois construction of Remark~\ref{rem:vallois}.) We assume for convenience that $x>0$. Then we consider the embedding constructed with $F(s) = s \wedge x$, and $h(s) = s$. The resulting embedding uses Az\'{e}ma-Yor to embed the left tail of the distribution, until the moment that the maximum first reaches $x$, and then uses the Vallois construction to embed the rest of the distribution. The choice of the Az\'{e}ma-Yor embedding will not affect the optimality of the construction, but fits in well with the construction described in this paper. We note, however, that the construction is the unique construction that maximises $\E \Psi(L_T^X)$ for a convex function $\Psi$, and that in this instance $L_T^X = \left(\sup_{t \le T} B_t\right) \wedge x + L_T^x$.

\begin{theorem}
Suppose $x \ge 0$. The stopping time $T$ defined in \eqref{eqn:Tdefn} with $F(s) = s \wedge x$, and $h(s) = s$ maximises $\E \Psi(L_T^x)$ for any convex function $\Psi$, over all embeddings $\tau$ of $\mu$ in $B_t$ such that $B_{t \wedge \tau}$ is UI.
\end{theorem}

\begin{proof}
By the optimality of the Vallois construction (see \citet{Vallois:92}), the problem can be reduced to a question of determinining which parts of $\mu$ we must embed before we hit $x$. For simplicity, we shall shift the problem to considering the local time at $0$, subject to having started at some point below $0$. This leaves us having to embed the part of $\mu$ lying above $0$ and some mass below $0$ dominated by $\mu$. We shall also suppose that $\mu$ has a density on $\R$ given by $f(y)$, however it is clear that this assumption can be removed by considering a limiting argument. Consequently, $f$ satisfies
\[
\int f(y) \, dy = 1, \int y f(y) \, dy <0.
\]

We shall prove the theorem using a Lagrangian argument. Suppose the optimal choice of embedding has a joint law for stopped distribution of $(B_T,L_T)$ given by $\rho(b,l)$. Then the problem of finding the optimal embedding can be restated as:
\begin{align}
&\mbox{Maximise} & \int_0^\infty \int_{-\infty}^\infty \Psi(l) \rho(b,l) \, db \,dl &&\nonumber\\
&\mbox{subject to} & \int_0^\infty \rho(b,l) \, dl &  =   f(b) & b>0 \label{eqn:lambda1}\\
& & \int_0^\infty \rho(b,l) \, dl &  \le   f(b) & b<0\label{eqn:lambda2}\\
& & \int_0^{\infty} b \rho(b,l) \, db &  =  \half \int_{-\infty}^\infty \int_l^\infty \rho(b,u) \, du \, dl &\label{eqn:Royn1}\\
& & \int_{-\infty}^0 b \rho(b,l) \, db &  =  -\half \int_{-\infty}^\infty \int_l^\infty \rho(b,u) \, du \, dl &\label{eqn:Royn2}\\
&&\rho(b,l) & \ge 0. &\nonumber
\end{align}
The necessity and sufficiency of the conditions \eqref{eqn:Royn1} and \eqref{eqn:Royn2} follow from Theorem~I.1 of \citet{Roynette:02}. We shall construct the dual problem, and show that, for a suitable choice of the Lagrangian multipliers and the conjectured $\rho(b,l)$, the primal and dual have the same value.

We introduce the multiplier $\lambda(b)$ corresponding to \eqref{eqn:lambda1} for $b>0$ and \eqref{eqn:lambda2} for $b<0$, and $\gamma(l),\eta(l)$ corresponding to \eqref{eqn:Royn1} and \eqref{eqn:Royn2} respectively. The resulting Lagrangian is:
\[ \begin{split}
\iint \rho(b,l) &\left[ \Psi(l) - \lambda(b) - \gamma(l) b \indic{b>0} - \eta(l) b \indic{b<0} + \half \int_0^l (\gamma(u) - \eta(u)) \, du\right] \, db \, dl\\
&{}+\int \lambda(b) f(b) \, db.\end{split}
\]
The corresponding dual problem is therefore to minimise $\int \lambda(b) f(b) \, db$ subject to the term in the square bracket being non-positive, and $\lambda(b) \ge 0$ for $b \le 0$ (as a consequence of the inequality in \eqref{eqn:lambda2}). For simplicity, we assume also that $\Psi(0) = \Psi'(0) = 0$; the second assumption is possible since $\E L_T^x = 2\E(B_T-x)_+$, which is fixed for any UI embedding, and $\Psi$ is convex.

The choice of the Lagrangian multipliers will depend on the optimal construction. As such, we are able to conjecture the solution: let $\delta<0$ satisfy: $\int_\delta^\infty b f(b) \, db = 0$. We then use the Vallois construction to embed $\mu$ restricted to $[\delta,\infty)$: let $\alpha(l) \ge 0$ decreasing, and $\beta(l) \le 0$ increasing be the Vallois boundaries, so that $T = \inf\{t \ge 0: B_t \not\in (\beta(L_t),\alpha(L_t))\}$ embeds $\mu$ restricted to $[\delta,\infty)$, when there is mass $\int_\delta^\infty f(b) \,db$ at $0$.

Then define
\begin{eqnarray*}
A(m) & = & \half \int_0^m \left[ \frac{1}{\alpha(k)} + \frac{1}{|\beta(k)|} \right] \, dk;\\
\gamma(l) & = & \int_0^l \frac{e^{A(u)}}{\alpha(u)} \int_0^u \Psi''(m) e^{-A(m)} \, dm\, du;\\
\eta(l) & = & \int_0^l \frac{e^{A(u)}}{\beta(u)} \int_0^u \Psi''(m) e^{-A(m)} \, dm\, du;\\
\lambda(b) & = & \begin{cases} \int_0^{\alpha^{-1}(b)} \frac{e^{A(u)}(\alpha(u) - b)}{\alpha(u)} \int_0^u \Psi''(m) e^{-A(m)} \, dm\, du &:b>0\\
\int_0^{\beta^{-1}(b)} \frac{e^{A(u)}(b-\beta(u))}{|\beta(u)|} \int_0^u \Psi''(m) e^{-A(m)} \, dm\, du &:b<0 \end{cases}.
\end{eqnarray*}
Note in particular that $\lambda(b) \ge 0$ for $b<0$. With these definitions, and using the properties of the functions $\alpha(l)$ and $\beta(l)$, it is now possible to show that
\[
\Psi(l) - \lambda(b) - \gamma(l) b \indic{b>0} - \eta(l) b \indic{b<0} + \half \int_0^l (\gamma(u) - \eta(u)) \, du \le 0,
\]
so that we have a feasible solution to the dual problem, and also
\[
\int \lambda(b) f(b) \, db = \E \Psi(L_T)
\]
when $T$ is the candidate embedding. Consequently the embedding is indeed optimal.
\end{proof}


%

\bibliography{general}

\begin{thebibliography}{20}
\expandafter\ifx\csname natexlab\endcsname\relax\def\natexlab#1{#1}\fi
\expandafter\ifx\csname url\endcsname\relax
  \def\url#1{\texttt{#1}}\fi
\expandafter\ifx\csname urlprefix\endcsname\relax\def\urlprefix{URL }\fi

\bibitem[{Az{\'e}ma et~al.(1980)Az{\'e}ma, Gundy, and Yor}]{AzemaGundyYor:79}
Az{\'e}ma, J., Gundy, R.~F., Yor, M., 1980. Sur l'int\'egrabilit\'e uniforme
  des martingales continues. In: Seminar on Probability, XIV (Paris, 1978/1979)
  (French). Vol. 784 of Lecture Notes in Math. Springer, Berlin, pp. 53--61.

\bibitem[{Az{\'e}ma and Yor(1979{\natexlab{a}})}]{AzemaYor:79}
Az{\'e}ma, J., Yor, M., 1979{\natexlab{a}}. \abc {U}ne solution simple au
  probl\`eme de {S}korokhod. In: S\'eminaire de Probabilit\'es, XIII (Univ.
  Strasbourg, Strasbourg, 1977/78). Vol. 721 of Lecture Notes in Math.
  Springer, Berlin, pp. 90--115.

\bibitem[{Az{\'e}ma and Yor(1979{\natexlab{b}})}]{AzemaYor:79b}
Az{\'e}ma, J., Yor, M., 1979{\natexlab{b}}. Le probl\`eme de {S}korokhod:
  compl\'ements \`a ``{U}ne solution simple au probl\`eme de {S}korokhod''. In:
  S\'eminaire de Probabilit\'es, XIII (Univ. Strasbourg, Strasbourg, 1977/78).
  Vol. 721 of Lecture Notes in Math. Springer, Berlin, pp. 625--633.

\bibitem[{Barlow et~al.(2000)Barlow, Burdzy, Kaspi, and
  Mandelbaum}]{Barlowetal:00}
Barlow, M., Burdzy, K., Kaspi, H., Mandelbaum, A., 2000. Variably skewed
  {B}rownian motion. Electron. Comm. Probab. 5, 57--66 (electronic).

\bibitem[{Bass(1983)}]{Bass:83}
Bass, R.~F., 1983. Skorokhod imbedding via stochastic integrals. In: Seminar on
  probability, XVII. Vol. 986 of Lecture Notes in Math. Springer, Berlin, pp.
  221--224.

\bibitem[{Chacon(1977)}]{Chacon:77}
Chacon, R.~V., 1977. Potential processes. Trans. Amer. Math. Soc. 226, 39--58.

\bibitem[{Chacon and Walsh(1976)}]{ChaconWalsh:76}
Chacon, R.~V., Walsh, J.~B., 1976. One-dimensional potential embedding. In:
  S\'eminaire de Probabilit\'es, X (Pr\`emiere partie, Univ. Strasbourg,
  Strasbourg, ann\'ee universitaire 1974/1975). Springer, Berlin, pp. 19--23.
  Lecture Notes in Math., Vol. 511.

\bibitem[{Cox(2005)}]{Cox:05}
Cox, A. M.~G., 2005. Extending {C}hacon-{W}alsh: Minimality and generalised
  starting distributions, preprint, {\ttfamily
  http://www.arxiv.org/abs/math.PR/0503535}.

\bibitem[{Cox and Hobson(2003)}]{CoxHobson:03}
Cox, A. M.~G., Hobson, D.~G., 2003. Skorokhod embeddings, minimality and
  non-centred target distributions, preprint, {\ttfamily
  http://www.arxiv.org/abs/math.PR/0310403}.

\bibitem[{Dubins(1968)}]{Dubins:68}
Dubins, L.~E., 1968. On a theorem of {S}korohod. Ann. Math. Statist. 39,
  2094--2097.

\bibitem[{Harrison and Shepp(1981)}]{HarrisonShepp:81}
Harrison, J.~M., Shepp, L.~A., 1981. On skew {B}rownian motion. Ann. Probab.
  9~(2), 309--313.

\bibitem[{Meilijson(1983)}]{Meilijson:83}
Meilijson, I., 1983. On the {A}z\'ema-{Y}or stopping time. In: Seminar on
  probability, XVII. Vol. 986 of Lecture Notes in Math. Springer, Berlin, pp.
  225--226.

\bibitem[{Monroe(1972)}]{Monroe:72}
Monroe, I., 1972. On embedding right continuous martingales in {B}rownian
  motion. Ann. Math. Statist. 43, 1293--1311.

\bibitem[{Ob{\l}\'{o}j(2004)}]{Obloj:04}
Ob{\l}\'{o}j, J., 2004. The maximality principle revisited: on certain optimal
  stopping problems, preprint {\ttfamily arxiv.org:math.PR/0405252}.

\bibitem[{Revuz and Yor(1999)}]{RevuzYor:99}
Revuz, D., Yor, M., 1999. Continuous martingales and {B}rownian motion, 3rd
  Edition. Vol. 293 of Grundlehren der Mathematischen Wissenschaften
  [Fundamental Principles of Mathematical Sciences]. Springer-Verlag, Berlin.

\bibitem[{Root(1969)}]{Root:69}
Root, D.~H., 1969. The existence of certain stopping times on {B}rownian
  motion. Ann. Math. Statist. 40, 715--718.

\bibitem[{Roynette et~al.(2002)Roynette, Vallois, and Yor}]{Roynette:02}
Roynette, B., Vallois, P., Yor, M., 2002. A solution to {S}korokhod's embedding
  for linear {B}rownian motion and its local time. Studia Sci. Math. Hungar.
  39~(1-2), 97--127.

\bibitem[{Skorokhod(1965)}]{Skorokhod:65}
Skorokhod, A.~V., 1965. Studies in the theory of random processes. Translated
  from the Russian by Scripta Technica, Inc. Addison-Wesley Publishing Co.,
  Inc., Reading, Mass.

\bibitem[{Vallois(1983)}]{Vallois:83}
Vallois, P., 1983. Le probl\`eme de {S}korokhod sur {${\bf R}$}: une approche
  avec le temps local. In: Seminar on probability, XVII. Vol. 986 of Lecture
  Notes in Math. Springer, Berlin, pp. 227--239.

\bibitem[{Vallois(1992)}]{Vallois:92}
Vallois, P., 1992. Quelques in\'egalit\'es avec le temps local en zero du
  mouvement brownien. Stochastic Process. Appl. 41~(1), 117--155.

\end{thebibliography}
\bibliographystyle{elsart-harv}
\end{document}